\documentclass[ijoo]{informs4}
\RequirePackage{tgtermes}
\RequirePackage{newtxtext}
\RequirePackage{newtxmath}
\RequirePackage{bm}
\RequirePackage{endnotes}

\OneAndAHalfSpacedXII 

\usepackage{algorithm}
\usepackage{algpseudocode}
\usepackage{tikz}
\usepackage{comment}

\usepackage{natbib}
 \bibpunct[, ]{(}{)}{,}{a}{}{,}%
 \def\bibfont{\small}%
\EquationsNumberedThrough    

\TheoremsNumberedThrough     
\ECRepeatTheorems  %

\MANUSCRIPTNO{IJOO-0001-2024.00}

\begin{document}


\RUNAUTHOR{Bonami et al.}

\RUNTITLE{Cutting Planes for Binarized Network Flow Problems}

\TITLE{Cutting Planes for Binarized Network Flow Problems}

\ARTICLEAUTHORS{%
\AUTHOR{Pierre Bonami}
\AFF{
GUROBI, \EMAIL{bonami@gurobi.com}}

\AUTHOR{Sanjeeb Dash}
\AFF{Foundations of Optimization \& Probability,
IBM, \EMAIL{sanjeebd@us.ibm.com}}

\AUTHOR{Anton Derkach}
\AFF{ORIE,
Cornell, \EMAIL{ad2265@cornell.edu}}

\AUTHOR{Andrea Lodi}
\AFF{Jacobs Technion-Cornell Institute,
Cornell Tech and Technion - IIT, \EMAIL{andrea.lodi@cornell.edu}}
} 

\ABSTRACT{

\noindent We consider integer programming problems with bounded general-integer variables belonging to the general class of network flow problems. For those, we computationally investigate the effect on mixed-integer linear programming (MIP) solvers of the different ways of producing extended formulations that replace a bounded general integer variable by a linear combination of a set of auxiliary binary variables linked by additional linear constraints. We show that MILP solvers perform very differently depending on which extended formulations is used and we interpret that different performance through the lens of cutting planes generation. Finally, we discuss a simple family of mixed-integer rounding inequalities that especially benefit from the reformulation, and we show its benefit within different MIP solvers. This provides methodological and practical guidelines for the use of those extended formulations in MIP and, to the best of our knowledge, this is the first extensive computational analysis of the topic. All our data and tables are publicly available at \url{https://github.com/anton-derkach1/binarizations}.}%




\KEYWORDS{Extended formulations, Integer Programming, Separation}

\maketitle


\section{Introduction}\label{sec:Intro}
Extended formulations (EF) form an important tool for mixed-integer programming (MIP). Informally, one builds an extended formulation of an MIP instance by adding auxiliary variables and linear constraints that link the original and auxiliary variables. See, e.g., \cite{U2011} for a discussion.
Even when the auxiliary variables are continuous and the linear programming (LP) relaxation of the extended formulation is no stronger than that of the original formulation (with respect to the optimal objective function value), it can happen that a family of cutting planes is stronger on the extended formulation than on the original one \citep{bodur}.
A well-known EF approach is to replace a bounded general integer variable by a linear combination of a set of auxiliary binary variables linked by additional linear constraints. There are several possible ways of making such a replacement.
For example, given an integer variable $x \in \{0,\ldots a\}$, one can perform what is referred to in the literature as \emph{full binarization} \citep{DGH18}, namely, add the constraints

\begin{equation}
\label{eq:fb}
    \sum_{k=1}^{a} k z^k = x,~
    \sum_{k=1}^{a} z^k   \leq 1, ~
    z^k \in \{0,1\},~k = 1, \ldots, a.
\end{equation}
Another type of binarization, called \emph{logarithmic binarization}, namely,
\begin{equation}
\label{eq:lb}
    \sum_{k=1}^{\lceil \log_2 (a+1)\rceil} 2^{k-1} z^k = x,~
    x \leq a,~
    z^k \in \{0,1\}, k = 1, \ldots, \lceil \log_2 (a+1) \rceil,
\end{equation}is popular in the literature as it has the smallest possible number of auxiliary binary variables.
\cite{OM2002} analyze both the full and logarithmic binarizations and show that pure 0-1 branching (without cut generation) would perform worse on the extended formulation than on the original MIP. Accordingly, they suggested that such extended formulations were unlikely to be useful in practice.

Several recent papers argue that binarization may enable the generation of more effective cutting planes to strengthen the LP relaxation of an MIP. In \cite{U2011}, full binarization is used as an intermediate step to define different general integer variables that, for fixed-charge network flow problems (and alike), link capacities, loads and times. 
Cuts on the resulting extended formulation are quite effective. \citet{BM2015} use a different binarization, namely
\begin{equation}
\label{eq:ub}
    \sum_{k=1}^{a} z^k = x,~
    1 \geq z^1 \geq z^2 \geq \cdots \geq z^a \geq 0,~
    z^k \in \{0,1\}, k = 1, \ldots, a,
\end{equation}
called \emph{unary binarization} in \cite{DGH18}, and analyze the effectiveness of cutting planes in this extended formulation both theoretically and computationally. \cite{DGH18} compare the quality of the above binarizations by comparing their split closure \citep{dash2010mir} and show that the unary and full binarizations are equivalent and dominate the logarithmic one with respect to the strength of split cuts that can be generated from these formulations. Finally, \cite{AVV2017} define an extended formulation of a class of fixed-charge transportation problems (similar to the problems already discussed in \cite{U2011}). They use full binarization strengthened by additional constraints linking the original and auxiliary variables. They show that a general-purpose MIP solver, namely  CPLEX, automatically exploits such a reformulation and solves a number of instances  
 orders of magnitude faster than with the original formulation (without binarization) mainly due to cutting plane generation.

\paragraph{Paper Contribution.}
In this paper, we aim to compare, from a computational standpoint, different binarizations and how well they can be exploited by state-of-the-art MIP solvers. In this sense, in Section \ref{sec:CPLEX4FCT}, we extend the results in \citep{AVV2017} and \citep{DGH18} by \emph{(i)} comparing CPLEX performance on multiple binarizations, i.e., the full, unary, and logarithmic binarizations, and \emph{(ii)} analyzing performance variations resulting from the way these extended formulations are expressed as inequality systems.
Moreover, in Section \ref{sec:form_cuts}, (iii) we show that a very simple and effective family of mixed-integer rounding (MIR) inequalities \cite{dash2010mir} can be separated in the MIPs discussed by \citep{U2011, AVV2017} by taking into account the flow-conservation constraints rewritten in terms of the binary $z$ variables, thus highlighting which inequalities play a central role in the effectiveness of the MIP solvers. 
Finally, in Section \ref{sec:MIP_solvers}, (iv) we show that these simple MIR cuts enhance CPLEX performance overall, and that these results are robust to the use of other solvers, namely Gurobi.
To the best of our knowledge, each one of the above contributions is new. In other words, this is the first extensive computational analysis of the topic and the results significantly contribute to clarify when and how binarizations can be effective in MIP with network flow structure, thus giving practical guidelines for their use. This is summarized in Section \ref{sec:conclusions}.

\section{Comparing different binarizations}
\label{sec:CPLEX4FCT}

In this section, we compare the impact of the different binarization strategies introduced in the previous section on the performance of a modern MIP solver, specifically CPLEX, while solving the fixed-charge transportation problem (FCT) in \cite{AVV2017}. The goal of this experiment is to use CPLEX in a black-box manner, i.e., just observe how the MIP solver performs on different formulations. In the next sections, our goal will be to interpret the results.

Given a bipartite graph $G=(V,E)$ with $V=N \cup M$, where $N$ and $M$ are suppliers and customers, respectively, the supply and demand capacities $s_i, d_j \in \mathbb{Z}_{+}$ for each supply and demand node $i \in N$ and $j \in M$, set the capacity of the arc $(i,j)\in E$ to $a_{ij}=\min\{s_i,d_j\}$. Given a cost vector $q \in \mathbb{Z}^{|N||M|}_{\geq 0}$, assume we want to minimize the total cost of using the arcs while satisfying the demand with the available supply. The resulting problem is an FCT problem:

\begin{subequations}
\label{eq:fct}
\begin{alignat}{2}
    (\text{FCT})\quad
    \text{min}\quad
        & q^T y \nonumber \\
    \text{s.t.}\quad
        & \sum_{j\in M} x_{ij} \le s_i
        && \quad i \in N \label{eq:fct1} \\
        & \sum_{i\in N} x_{ij} = d_j
        && \quad j \in M \label{eq:fct2} \\
        & y_{ij} \le x_{ij} \le a_{ij} y_{ij}
        && \quad (i,j) \in E \label{eq:fct3a} \\
        & y_{ij} \in \{0,1\}
        && \quad (i,j) \in E. \label{eq:fct3}
\end{alignat}
\end{subequations}
Note that, for fixed $y \in \{0,1\}^{|E|}$, the constraint matrix on $x$ is totally unimodular, and thus constraints \eqref{eq:fct3} and the integrality of $a_{ij}$ imply that $x$ is integral.

\cite{AVV2017} define an extended formulation of the FCT formulation above as follows. They apply full binarization \eqref{eq:fb}  to the $x_{ij}$ variables in (\ref{eq:fct}) and add additional valid constraints, namely,
\begin{equation}
\label{eq:avv_s}
    \sum_{k=1}^{a_{ij}} z_{ij}^k = y_{ij} \quad (i,j) \in E, 
\end{equation}
which are based on constraints (\ref{eq:fct3a}) relating the flow variables $x_{ij}$, the binary variables $y_{ij}$ indicating whether the arc $(i,j)$ is used, and the capacity $a_{ij}$ of the arc $(i,j)$. The resulting formulation is:
\begin{subequations}
\label{eq:avv}
\begin{alignat}{3}
    (\text{AvV})\quad
    \text{min}\quad
        & q^T y \nonumber \\
    \text{s.t.}\quad
        & \sum_{j\in M} \sum_{k=0}^{a_{ij}} k\, z_{ij}^k \le s_i
        && \quad i \in N \label{eq:avv1} \\
        & \sum_{i\in N} \sum_{k=0}^{a_{ij}} k\, z_{ij}^k = d_j
        && \quad j \in M \label{eq:avv2} \\
        & y_{ij} \le x_{ij} \le a_{ij} y_{ij}
        && \quad (i,j) \in E \label{eq:avv3} \\
        & \sum_{k=0}^{a_{ij}} k\, z_{ij}^k = x_{ij}
        && \quad (i,j) \in E \label{eq:avv4} \\
        & \sum_{k=1}^{a_{ij}} z_{ij}^k = y_{ij}
        && \quad (i,j) \in E \label{eq:avv5} \\
        & \sum_{k=0}^{a_{ij}} z_{ij}^k = 1
        && \quad (i,j) \in E \label{eq:avv6} \\
        & y_{ij} \in \{0,1\}
        && \quad (i,j) \in E \label{eq:avv7} \\
        & z_{ij}^k \in \{0,1\}
        && \quad (i,j) \in E,\ 0 \le k \le a_{ij}. \label{eq:avv8}
\end{alignat}
\end{subequations}
Notice that (\ref{eq:avv1}) and (\ref{eq:avv2}) are just the flow-conservation constraints (\ref{eq:fct1}) and (\ref{eq:fct2}) rewritten in terms of $z$. \cite{AVV2017} observe that even though the projection of the linear relaxation of (AvV) onto $(x,y)$ is equal to the linear relaxation of (FCT), empirically, MIP solvers perform much better on (AvV) 
than on (FCT). Furthermore, we will show that constraints (\ref{eq:avv_s}) strengthen the full binarization.  

Table \ref{tab:default_cplex} reports results obtained by solving 8 different formulations with CPLEX: 
\begin{enumerate}
    \item \textbf{FCT}: the FCT original formulation in \cite{AVV2017}, with no binarization.
    \item \textbf{FullB}: the full binarization \eqref{eq:fb}.
    \item \textbf{AvV}: the full binarization plus the strengthening of constraints \eqref{eq:avv_s}.
    \item \textbf{UnaryB$^+$}: the unary binarization \eqref{eq:ub} with the strengthened constraints (similar to (\eqref{eq:avv_s}):
    \begin{equation}
    \label{eq:ub_s}
    z_{ij}^1 = y_{ij} \quad (i,j) \in E.
    \end{equation}
    \item \textbf{LogB$^+$}: the logarithmic binarization \eqref{eq:lb} with the strengthened constraints (similar to (\eqref{eq:avv_s}):
    \begin{equation}
    \label{eq:lb_s}
   \sum_{k=1}^{\lceil \log_2 (a_{ij}+1)\rceil} z_{ij}^k \geq y_{ij} \quad (i,j) \in E.
    \end{equation}
    \item \textbf{AvV - z}: the \textbf{AvV} formulation in which -- given our hypothesis that the flow-conservation constraints (\ref{eq:avv1}) and (\ref{eq:avv2}) are a rich source of cuts responsible for better empirical performance of (AvV) compared to (FCT) (see \cite{AVV2017}) -- we remove constraints (\ref{eq:avv1})-(\ref{eq:avv2}).
    \item \textbf{AvV + U}: the \textbf{AvV} formulation amended by variables and constraints suggested by \cite{U2011}, where after the FCT problem is fully binarized, new general integer variables are introduced to link equally sized flows from all supply (resp. demand) nodes into a single demand (resp. supply) node. More precisely, set $C=\max\{a_{ij}: (i,j)\in E\}$ and define aggregated integer variables $u_i^k$ and $w_j^k$ as follows:
     \begin{subequations}
    \label{eq:uchoa_agg}
    \begin{align}
   u_i^k & = \sum_{j: k \leq a_{ij}}z_{ij}^k \quad  k =1,\dots, C \label{eq:uchoa_agg1} \\
   w_j^k & = \sum_{i: k \leq a_{ij}}z_{ij}^k \quad k =1,\dots, C \label{eq:uchoa_agg2}.
    \end{align}
    \end{subequations}
    The flow-conservation constraints (\ref{eq:fct1})-(\ref{eq:fct2}) can be written in terms of  aggregated variables as
    \begin{subequations}
    \label{eq:uchoa_flows}
    \begin{align}
   \sum_{k=1}^C k  \, u_i^k \leq s_i& \quad i \in N  \label{eq:uchoa_flows1} \\
  \sum_{k=1}^C k  \, w_j^k = d_j & \quad j \in M \label{eq:uchoa_flows2}.
    \end{align}
    \end{subequations}
    In \cite{U2011}, it is shown that constraints (\ref{eq:uchoa_flows}) are a good source for cuts. So, (AvV + U) is defined as (AvV) with (\ref{eq:uchoa_agg}) and (\ref{eq:uchoa_flows}).
    \item \textbf{AvV + U - z}: the (AvV + U) formulation, where, to isolate the effect of flow-conservation constraints in terms of aggregated variables  (\ref{eq:uchoa_agg}), we remove the flow-conservation constraints in terms of binary $z$ variables (\ref{eq:avv1}) and (\ref{eq:avv2}). Hence, this formulation contains only flow-conservation constraints in terms of aggregated variables (\ref{eq:uchoa_agg}), i.e., constraints (\ref{eq:uchoa_flows}).
\end{enumerate}

In this computational experiment, we use a dataset of 20 randomly generated instances. The number of suppliers and customers is the same and equal to $n$, where $n \in \{30,40\}$. Capacities are generated with a parameter $B\in \{10,20\}$ using the procedure described in \cite{AVV2017}, where $B$ serves as the ceiling for all capacities. We use the same factor $r=0.95$ as the total demand to total supply ratio. We generate 5 instances for each combination of $n,B$, and $r$.
Table~\ref{tab:default_cplex} reports the results of our experiments using the following metrics: (1) the instance classes $(n,\,B)$; (2) LP gap\%, which is computed as 
$100 \times (\text{OPT}-\text{LP relaxation value})/\text{OPT}$; (3) Root gap\%, which is computed as
$100 \times (\text{OPT}-\text{Root value})/\text{OPT}$; (4) the number of branch-and-bound nodes explored; and (5) the solution time in CPU seconds. 
For each class $(n, \, B)$, the results are averaged over the 5 instances, and the Avg column contains the averages of the values of each parameter class $(n, \, B)$. 
We perform all experiments with CPLEX using version 22.1.2 on a machine with 128 GBytes of memory,
and four 2.8 GHz Intel Xeon E7-4890 v2 processors, each with 15 cores, for a total of 60 cores. However, all solver runs use 8 threads with a time limit of 3600 seconds.


\setlength{\tabcolsep}{4pt} 
\renewcommand{\arraystretch}{0.65} 

\begin{table}[htbp]
\centering
\caption{\textbf{FCT formulations compared using CPLEX.}\label{tab:default_cplex}}

{\footnotesize

\begin{tabular}{lr|rrrrr}
\hline
& \up\down $(n,\,B)$ & (30,\,10) & (30,\,20) & (40,\,10) & (40,\,20) & Avg \\
\hline
& LP gap\% & 13.76 & 14.29 & 12.00 & 15.01 & 13.77 \\
\hline
\textbf{FCT} & 
\quad Root Gap\% & 2.42 & 3.39 & 2.15 & 4.51 & 3.12 \\
& \quad Nodes      & 8,980.2 & 21,348.2 & 14,684.2 & 518,699.0 & 140,927.9 \\
& \quad Time       & 7.13 & 19.10 & 21.64 & 1,689.86 & 434.43 \\
\hline
\textbf{FullB} & 
\quad Root Gap\% & 0.00 & 0.10 & 0.00 & 0.22 & 0.08 \\
& \quad Nodes      & 0.0 & 14.4 & 0.0 & 386.6 & 100.3 \\
& \quad Time       & 1.65 & 4.09 & 3.96 & 12.91 & 5.65 \\
\hline
\textbf{AvV} & 
\quad Root Gap\% & 0.00 & 0.04 & 0.00 & 0.11 & 0.04 \\
& \quad Nodes      & 0.0 & 3.8 & 0.0 & 22.4 & 6.6 \\
& \quad Time       & 0.66 & 1.88 & 1.33 & 4.88 & 2.19 \\
\hline
\textbf{UnaryB$^+$} & 
\quad Root Gap\% & 7.70 & 10.84 & 8.29 & 12.36 & 9.80 \\
& \quad Nodes      & 964,661.8 & 1,818,001.0 & 1,418,344.6 & 931,969.2 & 1,283,244.2 \\
& \quad Time       & 1,574.17 & 3,601.33 & 3,162.92 & 3,600.52 & 2,984.73 \\
\hline
\textbf{LogB$^+$} & 
\quad Root Gap\% & 5.58 & 8.17 & 4.90 & 9.57 & 7.06 \\
& \quad Nodes      & 303,996.0 & 1,623,406.6 & 651,503.8 & 1,022,860.6 & 900,441.8 \\
& \quad Time       & 379.38 & 3,364.45 & 1,607.92 & 3,604.00 & 2,238.93 \\
\hline
\textbf{AvV - z} & 
\quad Root Gap\% & 2.13 & 3.71 & 2.15 & 4.63 & 3.15 \\
& \quad Nodes      & 3,378.0 & 5,744.0 & 4,830.2 & 67,597.2 & 20,387.4 \\
& \quad Time       & 9.81 & 33.17 & 22.18 & 523.22 & 147.10 \\
\hline
\textbf{AvV + U} & 
\quad Root Gap\% & 0.00 & 0.00 & 0.00 & 0.20 & 0.05 \\
& \quad Nodes      & 0.0 & 0.0 & 0.0 & 75.8 & 19.0 \\
& \quad Time       & 0.73 & 1.97 & 1.27 & 4.61 & 2.14 \\
\hline
\textbf{AvV + U - z} & 
\quad Root Gap\% & 0.00 & 0.00 & 0.00 & 0.17 & 0.04 \\
& \quad Nodes      & 0.0 & 0.0 & 0.0 & 79.6 & 19.9 \\
& \quad Time       & 0.58 & 1.98 & 1.43 & 5.24 & 2.31 \\
\hline
\end{tabular}
}
\end{table}

\vspace{-0.5cm}
Table~\ref{tab:default_cplex} demonstrates significant performance differences in CPLEX among the analyzed formulations despite having the same LP relaxations (shown in \cite{AVV2017}). The full binarization significantly outperforms (FCT): it reduces the root gap\% by two orders of magnitude (on average, 0.08\% vs 3.12\%) and, for some instances, even closes the gap solving the problem at the root node. Consequently, the number of explored nodes is reduced by at least a factor of $10^3$, greatly decreasing the computing time. 

The difference in performance between (AvV) and the full binarization due to strengthening (\ref{eq:avv_s}) is particularly interesting. We see that the strengthening, which does not affect the LP bound, reduces the number of explored nodes for class 
(40,\,20) by more than one order of magnitude, as well as resulting in uniformly better root gap\%, number of explored nodes, and solution times for all other parameter classes.

Even after strengthening, the unary and logarithmic formulations experience a severe degradation in performance relative to (FCT) in root gap\%, number of explored nodes, and solution times, hitting the time limit of 3,600 seconds for many instances. The difference in CPLEX performance on the full and unary binarizations is striking given their theoretical equivalence with respect to split closures \citep{DGH18} and the fact that many of the cuts generated CPLEX are split cuts. Our results indicate that CPLEX finds strong cutting planes for the full formulation but cannot find the equivalent cuts for the unary formulation. One of the main algorithms for separating split cuts in modern MIP solvers is the MIR algorithm of \cite{MarchandWolsey2001}. We will show in the next section that strong split cuts can be derived by applying the MIR formula directly to specific rows of the full formulation. By contrast, to get an equivalent cut in the unary formulation one would first need to perform a specific aggregation on which to apply an MIR formula. In algorithms such as \cite{MarchandWolsey2001}, the construction of such an aggregation is heuristic, and this is why CPLEX is not able to construct these cuts.

The severe degradation in performance of (AvV - z) compared to (AvV) shows that removing flow constraints in terms of binary variables $z$ from (AvV) and the full binarization prevents CPLEX from exploiting the reformulation effectively. As we will show in the next section, effective cutting planes can be separated from individual flow-conservation constraints (\ref{eq:avv1}) and (\ref{eq:avv2}), so removing those constraints heavily degrades CPLEX performance because it is unable to perform the appropriate aggregations to obtain the flow-conservation constraints written in terms of the $z$ variables (even if such aggregations exist \citep{dash2010mir}).

The performances of (AvV + U) and (AvV + U - z) are comparable to that of (AvV) with the exception of the instances in the class (40,\,20), where (AvV + U) and (AvV + U - z) result in 3 times more nodes than for (AvV). Given the downside of including additional variables and constraints in the formulation (i.e., larger formulations), the empirical results imply no benefit of formulating the FCT problem as (AvV + U) or (AvV + U - z) over (AvV).

\section{MIR cuts from flow-conservation constraints}
\label{sec:form_cuts}

In this section, we show that a simple and effective family of cuts can be separated on certain binarizations of the fixed-charge transportation problem \cite{AVV2017} described in the previous section and of another class of MIPs, namely the capacitated minimum spanning tree problem  \cite{U2011}.

\paragraph{Capacitated Minimum Spanning Tree (CMST).} Let $G=(V,A)$ be a directed graph with vertices $V=\{0,1,\dots, n\}$ and $m=|A|$ arcs. Vertex 0 is the \emph{root}. Each non-root vertex $i$ has a positive demand $d_i \in \mathbb{Z}$. When all such demands are equal to 1, the instance is said to be unitary. Vertex 0 has a demand $d_0$ equal to 0. Each arc $a \in A$ has a nonnegative cost $q_a$. Given a positive integer $C$ greater than or equal to the maximum demand, the CMST objective is to find a minimum cost spanning tree for $G$ such that the total demand of the vertices in each subtree hanging from the root does not exceed $C$. Let binary variables $z_a^k$ indicate whether arc $a=(i,j)$ belongs to the optimal arborescence and that the total demand of all vertices in the sub-arborescence rooted in $j$ is exactly $k$. In other words, $z_a^k =1$ indicates that a flow of $k$ travels over arc $a$ and that a link of capacity $k$ needs to be installed over arc $a$. Note that variables $z_{ij}^k$ with $k>C-d_i$ do not have to be used. We define a directed multigraph $G_C=(V,A_C)$, where $A_C$ contains arcs $(i,j)^k$ for each $(i,j)\in A$ and $k=1,\dots, C-d_i$. Let $V_{+}=\{1,\dots, n\}$ be the set of non-root vertices and $\delta^{-}(i)$ and $\delta^{+}(i)$ be the sets of multigraph arcs going into vertex $i$ and coming out of vertex $i$, respectively. Analogously to the FCT problem, we write the AvV formulation as
\begin{subequations}
\label{eq:avv_cmst}
\begin{alignat}{3}
    (\text{AvV})\quad
    \text{min}\quad
        & \sum_{a^k \in A_C} q_a^k z_a^k \nonumber \\[2mm]
    \text{s.t.}\quad
        & \sum_{a^k \in \delta^{-}(i)} z_{a}^k = 1
        && \quad i \in V_{+} \label{eq:avv_cmst1} \\
        & \sum_{a^k \in \delta^{-}(i)} k\, z_{a}^k - \sum_{a^k \in \delta^{+}(i)} k\, z_{a}^k = d_i
        && \quad i \in V_{+} \label{eq:avv_cmst2} \\
        & z_{a}^k \in \{0,1\}
        && \quad a^k \in A. \label{eq:avv_cmst3}
\end{alignat}
\end{subequations}

Similarly to (AvV) for the FCT problem, (AvV) for the CMST problem contains flow-conservation constraints (\ref{eq:avv_cmst2}). The strengthening similar to (\ref{eq:avv5}) is implicit in the objective function. Additionally, the formulation has in-degree constraints (\ref{eq:avv_cmst1}). Analogously to the FCT problem, we can construct formulations (AvV + U) and (AvV + U - z) for CMST using general integer variables to aggregate inflows and outflows of the same size. In the computational experiments performed on the CMST problem, we use a dataset of 20 CMST instances from the ORLIB library. The tc80 and te80 instances have $n=80$ non-root vertices with unitary demands and capacity $C=5$.

\paragraph{Formulation cuts.} In this section, we investigate the effectiveness of cuts based on the rows of formulation, also known as \emph{formulation cuts} in \cite{dash2010twostep}. Given a set
\begin{equation}
    P = \left\{v \in \mathbb{R}^{|J|}, \,x\in \mathbb{Z}^{|I|}: \sum_{j \in J} c_jv_j + \sum_{i \in I} a_ix_i \geq b, \, v, \, x \geq 0 \right\} \nonumber
\end{equation}
with the {\em base} inequality $\sum_{j \in J} c_jv_j + \sum_{i \in I} a_ix_i \geq b$, the MIR inequality
\begin{equation}
    \sum_{j\in J} \max\{c_j, 0\}v_j + \sum_{i \in I} (\hat b\lfloor a_i \rfloor + \min\{\hat b, \hat a_i\})x_i \geq \hat b \lceil b \rceil \nonumber
\end{equation}
is valid for $P$, where $\hat b = b - \lfloor b \rfloor$ and $\hat a_i = a_i - \lfloor a_i \rfloor$.
For a given row of an MIP formulation and a given fractional solution, the procedure in \cite{dash2010twostep} generates base inequalities by dividing the row by the coefficient of an integer variable that is currently fractional. Variables with upper bounds are typically complemented if their current value is closer to their upper bound than their lower bound. After writing the MIR inequality, the complementing is undone to obtain the cut in the original space. This procedure produces effective MIR cuts. In addition to this procedure, we also generate formulation cuts where we do not complement any variables. 

The goal of the experiments in this section is to show the effectiveness of formulation cuts for both FCT and CMST instances where the the full binarization plus strengthening has been applied, i.e., the (AvV) formulation. In order to do that, we apply a classical cutting plane approach in which we separate rank-1 formulation cuts in rounds (see, e.g., \cite{dash2010mir}) till no more formulation cuts can be separated. No other cuts are used. To measure the effect of formulation cuts, we collect the final LP gap\% after adding cuts and the number of added cuts for each instance $(n,\,B,\,s)$. We compare the LP gap\% after adding cuts for each formulation to the initial LP gap\% and root gap\% of the AvV formulation using CPLEX default cut generation settings are used as benchmarks. We perform the experiment on all the formulations previously considered except the (FCT) original formulation. In addition, we also test the effect of formulation cuts on the full binarization where the first round of cuts separated consists of Gomory Mixed-integer (GMI) cuts \citep{gomory1960rm}.
We report the results for FCT and CMST in Tables \ref{tab:fct_formulation} and \ref{tab:cmst_formulation}, respectively.


\setlength{\tabcolsep}{2pt}
\renewcommand{\arraystretch}{1.0}

\begin{table}[htbp]
\footnotesize
\centering
\caption{FCT formulation cut effectiveness comparisons.}
\label{tab:fct_formulation}
\resizebox{\textwidth}{!}{%
\begin{tabular}{lr r|
                rr|rr|rr|rr|rr|rr|rr|rr}
\hline
\up\down $(n,\,B,\,s)$ & {\begin{tabular}{@{}c@{}}LP\\ gap\% \end{tabular}} & 
{\begin{tabular}{@{}c@{}}Root \\gap\% \end{tabular}} &
\multicolumn{2}{c|}{FullB} &
\multicolumn{2}{c|}{\begin{tabular}{@{}c@{}}FullB with \\ GMI cuts\end{tabular}} &
\multicolumn{2}{c|}{AvV} &
\multicolumn{2}{c|}{AvV - z} &
\multicolumn{2}{c|}{UnaryB$^+$}  &
\multicolumn{2}{c|}{LogB$^+$} &
\multicolumn{2}{c|}{\begin{tabular}{@{}c@{}}AvV + U\end{tabular}} &
\multicolumn{2}{c}{\begin{tabular}{@{}c@{}}AvV + U - z\end{tabular}}\\
\cline{4-5}\cline{6-7}\cline{8-9}\cline{10-11}\cline{12-13}\cline{14-15}\cline{16-17}\cline{18-19}
 & & &
\begin{tabular}{@{}c@{}}LP Gap\%\\w. cuts\end{tabular} & \begin{tabular}{@{}c@{}}Cuts\\Added\end{tabular} &
\begin{tabular}{@{}c@{}}LP Gap\%\\w. cuts\end{tabular} & \begin{tabular}{@{}c@{}}Cuts\\Added\end{tabular} &
\begin{tabular}{@{}c@{}}LP Gap\%\\w. cuts\end{tabular} & \begin{tabular}{@{}c@{}}Cuts\\Added\end{tabular} &
\begin{tabular}{@{}c@{}}LP Gap\%\\w. cuts\end{tabular} & \begin{tabular}{@{}c@{}}Cuts\\Added\end{tabular} &
\begin{tabular}{@{}c@{}}LP Gap\%\\w. cuts\end{tabular} & \begin{tabular}{@{}c@{}}Cuts\\Added\end{tabular} &
\begin{tabular}{@{}c@{}}LP Gap\%\\w. cuts\end{tabular} & \begin{tabular}{@{}c@{}}Cuts\\Added\end{tabular} &
\begin{tabular}{@{}c@{}}LP Gap\%\\w. cuts\end{tabular} & \begin{tabular}{@{}c@{}}Cuts\\Added\end{tabular} &
\begin{tabular}{@{}c@{}}LP Gap\%\\w. cuts\end{tabular} & \begin{tabular}{@{}c@{}}Cuts\\Added\end{tabular} \\
\hline
$(30,\, 10,\, 1)$ & 13.73 & 0.00 & 13.73 & 70 &  12.51 & 150 & 0.56 & 142
& 13.73 & 0 & 13.73 & 0 & 13.73 & 72 & 0.56 & 250 & 3.16&87 \\ 
$(30,\, 10,\, 2)$ & 14.35 & 0.00 & 14.35 & 129 &  13.77 & 159 & 1.22 & 124
& 14.35 & 0 & 14.35 & 0 & 14.35 & 107 & 1.00 & 246 & 2.51& 102\\ 
$(30,\, 10,\, 3)$ & 15.79 & 0.00 & 15.79 & 93 &  15.15 & 170 & 1.35 & 162
& 15.79 & 0 & 15.79 & 0 & 15.79 & 74 & 1.36 & 281 &3.94 & 83\\ 
$(30,\, 10,\, 4)$ & 12.35 & 0.00 & 12.35 & 79 &  11.31 & 178 & 1.11 & 85
& 12.35 & 0 & 12.35 & 0 & 12.35 & 70 & 1.11 & 163 & 1.56& 81\\ 
$(30,\, 10,\, 5)$ & 12.61 & 0.00 & 12.61 & 91 &  12.09 & 148 & 0.90 & 141
& 12.61 & 0 & 12.61 & 0 & 12.61 & 82 & 0.90 & 226 &3.36 & 74\\ \hline
$(30,\, 20,\, 1)$ & 15.77 & 0.00 & 15.77 & 178 &  15.08 & 349 & 2.69 & 227
& 15.77 & 0 & 15.77 & 0 & 15.77 & 155 & 2.69 & 414 & 4.40& 139\\ 
$(30,\, 20,\, 2)$ & 13.41 & 0.00 & 13.41 & 132 & 12.46 & 229 &  1.21 & 177
& 13.41 & 0 & 13.41 & 0 & 13.41 & 81 & 1.21 & 351 &2.41 &170 \\ 
$(30,\, 20,\, 3)$ & 14.04 & 0.00 & 14.04 & 303 &  13.64 & 302 & 0.69 & 211
& 14.04 & 0 & 14.04 & 0 & 14.04 & 155 & 0.56 & 410 & 2.23& 149\\ 
$(30,\, 20,\, 4)$ & 13.51 & 0.00 & 13.51 & 154 &  12.70 & 280 & 2.49 & 166
& 13.51 & 0 & 13.51 & 0 & 13.51 & 166 & 2.60 & 325 & 3.59& 121\\ 
$(30,\, 20,\, 5)$ & 14.72 & 0.00 & 14.72 & 154 & 13.60 & 262 &  1.12 & 207
& 14.72 & 0 & 14.72 & 0 & 14.72 & 158 & 1.07 & 373 & 2.48& 116\\ \hline
$(40,\, 10,\, 1)$ & 12.62 & 0.00 & 12.62 & 115 & 11.63 & 231 &  0.52 & 184
& 12.62 & 0 & 12.62 & 0 & 12.62 & 107 & 0.50 & 323 &2.84 & 117\\ 
$(40,\, 10,\, 2)$ & 14.19 & 0.01 & 14.19 & 108 & 13.44 & 193 &  2.02 & 177
& 14.19 & 0 & 14.19 & 0 & 14.19 & 99 & 2.03 & 328 & 4.43& 123\\ 
$(40,\, 10,\, 3)$ & 11.63 & 0.00 & 11.63 & 130 &  10.95 & 172 & 0.81 & 116
& 11.63 & 0 & 11.63 & 0 & 11.63 & 148 & 0.81 & 214 & 2.37& 90\\ 
$(40,\, 10,\, 4)$ & 10.32 & 0.00 & 10.32 & 113 &  9.40 & 263 & 0.20 & 140
& 10.32 & 0 & 10.32 & 0 & 10.32 & 106 & 0.21 & 253 & 0.68&101 \\ 
$(40,\, 10,\, 5)$ & 11.22 & 0.01 & 11.22 & 97 &  10.81 & 209 & 0.97 & 130
& 11.22 & 0 & 11.22 & 0 & 11.22 & 92 & 0.97 & 244 & 2.21& 100\\ \hline
$(40,\, 20,\, 1)$ & 14.62 & 0.00 & 14.62 & 258 &  14.22 & 310 & 1.22 & 273
& 14.62 & 0 & 14.62 & 0 & 14.62 & 255 & 1.43 & 508 & 2.80&184 \\ 
$(40,\, 20,\, 2)$ & 16.59 & 0.40 & 16.59 & 207 &  16.24 & 356 & 1.55 & 352
& 16.59 & 0 & 16.59 & 0 & 16.59 & 236 & 1.55 & 629 & 3.61&237 \\ 
$(40,\, 20,\, 3)$ & 16.45 & 4.19 & 16.45 & 193 &  15.66 & 276 & 2.83 & 263
& 16.45 & 0 & 16.45 & 0 & 16.45 & 143 & 2.99 & 490 & 4.92&187 \\ $(40,\, 20,\, 4)$ & 12.04 & 0.01 & 12.04 & 187 &  11.53 & 386 & 1.48 & 267
& 12.04 & 0 & 12.04 & 0 & 12.04 & 151 & 1.53 & 498 & 3.09&168 \\ 
$(40,\, 20,\, 5)$ & 15.35 & 0.27 & 15.35 & 227 &  14.83 & 368 & 1.45 & 278
& 15.35 & 0 & 15.35 & 0 & 15.35 & 208 & 1.45 & 524 & 2.66& 198 \\ \hline
\hline
\textbf{Average} & 13.77 & 0.24 & 13.77 & 150.9 & 13.05 & 249.55 &  1.32 & 191.1 & 13.77 & 0 & 13.77 & 0 & 13.77 & 133.25 & 1.33 & 352.5 &2.96 & 131.35\\ \hline
\end{tabular}
}
\end{table}

The results in Table~\ref{tab:fct_formulation} show that formulation cuts are separated but do not move the bound in the full binarization, whereas they are extremely effective in reducing the gap from an average of 13.77\% to an average of 1.32\% for (AvV). This is not as good as the root gap\% obtained by CPLEX adding the entire arsenal of its cuts at the root node (0.24\%), but it is still remarkably close. Furthermore, adding one round of GMI cuts to the full binarization changes the bound only slightly leaving the binarization unable to exploit formulation cuts in the same capacity as (AvV). This difference in performance emphasizes the importance of the strengthening (\ref{eq:avv_s}) for effective separation of formulation cuts.

The results for (AvV - z) once again demonstrate the significance of rewriting flow constraints with $z$ variables in the ability of (formulation) cuts to strengthen the binarization. Essentially, simple MIR cuts from the formulation rows (and specifically from flow-conservation constraints) are effective (as somehow anticipated in the previous section) but are not separated by the CPLEX default separation procedures if (\ref{eq:avv1}) and (\ref{eq:avv2}) are not present.

Moreover, it is remarkable that formulation cuts were completely ineffective for the unary and logarithmic binarizations with strengthening leaving the gap unchanged.

Finally, we see that there are no significant differences in the gap closed after applying formulation cuts between (AvV) and (AvV + U). The drawback to using (AvV + U) is that, on average, a larger number of formulation cuts is added in (AvV + U) than (AvV) to obtain a similar average gap closed. We should note the difference between the average gaps of 2.96\% and 1.33\% after adding cuts for (AvV + U - z) and (AvV + U), respectively. This performance difference confirms the superiority of formulations containing flow-conservation constraints written in terms of binary variables $z$ for effective cut separation.

\begin{table}[htbp]
\centering
\footnotesize
\caption{CMST formulation cut effectiveness comparisons.}
\label{tab:cmst_formulation}
\resizebox{0.65\textwidth}{!}{ 
\begin{tabular}{lr r|
                rr|rr|rr}
\hline
Instance & {\begin{tabular}{@{}c@{}}LP\\ gap\% \end{tabular}} & 
{\begin{tabular}{@{}c@{}}Root \\gap\% \end{tabular}}  &
  \multicolumn{2}{c|}{AvV} &
  \multicolumn{2}{c|}{AvV + U} &
  \multicolumn{2}{c}{AvV + U - z} \\
\cline{4-5}\cline{6-7}\cline{8-9}
 & & &
\begin{tabular}{@{}c@{}}LP Gap\%\\w. cuts\end{tabular} & \begin{tabular}{@{}c@{}}Cuts\\Added\end{tabular} &
\begin{tabular}{@{}c@{}}LP Gap\%\\w. cuts\end{tabular} & \begin{tabular}{@{}c@{}}Cuts\\Added\end{tabular} &
\begin{tabular}{@{}c@{}}LP Gap\%\\w. cuts\end{tabular} & \begin{tabular}{@{}c@{}}Cuts\\Added\end{tabular} \\
\hline
tc80-1 & 6.82 & 0.79 & 0.53 & 319 & 0.56 & 616 & 0.73 & 275 \\ 
tc80-2 & 5.78 & 1.21 & 0.53 & 339 & 0.53 & 643 & 0.54 & 301  \\ 
tc80-3 & 5.82 & 0.85 & 0.75 & 304 & 0.75 & 607 &0.97 & 298\\ 
tc80-4 & 5.56 & 0.78 & 0.41 & 252 & 0.41 & 520 &0.44 & 382\\ 
tc80-5 & 6.45 & 0.80 & 0.45 & 372 & 0.47 & 663 & 0.59& 320\\ \hline
te80-1 & 6.01 & 1.07 & 1.71 & 248 & 1.70 & 489 & 1.95& 214\\ 
te80-2 & 6.92 & 1.17 & 2.06 & 290 & 2.06 & 535 & 2.40& 223\\ 
te80-3 & 6.13 & 1.05 & 1.64 & 272 & 1.61 & 522 & 1.81& 256\\ 
te80-4 & 7.35 & 1.64 & 1.51 & 330 & 1.51 & 653 &1.75 & 302\\ 
te80-5 & 7.11 & 1.69 & 2.40 & 315 & 2.39 & 585 & 2.75& 222\\ \hline \hline
\textbf{Average} & 6.40 & 1.11 & 1.20 & 304.1 & 1.20 & 583.3 &1.39 &269.3 \\ \hline
\end{tabular}
} 
\end{table}

Since the results in Table~\ref{tab:fct_formulation} show that formulation cuts are ineffective at reducing the gap for the full binarization with or without GMI cuts, (AvV - z), and the unary and logarithmic formulations with strengthening, we excluded them from Table~\ref{tab:cmst_formulation}. The results for CMST and FCT are similar. Formulation cuts close the gap from an average of 6.40\% to an average of 1.20\% for (AvV), comparable to the root gap\% average of 1.11\% that CPLEX obtains using its entire arsenal of cuts. Both (AvV) and (AvV + U) have the same average gap of 1.20\%, reinforcing the idea that there is no compelling reason to use (AvV + U) instead of (AvV) in this context. The slightly worse gap closed for (AvV + U - z) than for (AvV + U) supports the superiority of formulations containing flow-conservation constraints written in terms of binary variables $z$ for effective separation of formulation cuts.

\section{MIP solvers and formulation cuts}
\label{sec:MIP_solvers}

The goal of this section is twofold. On the one side, we show that both FCT and CMST, (AvV) strengthened by formulation cuts separated as described in the previous section, can be solved very effectively by CPLEX and that formulation cuts account for most of the CPLEX cut power as observed by \cite{AVV2017}. On the other side, we show that our findings are robust to the use of different MIP solvers by repeating the same experiments with Gurobi. 

To achieve those goals, we solve for both FCT (Table \ref{tab:avv_mir_form}) and CMST (Table \ref{tab:avv_zerohalf_form}), the following formulations using CPLEX:
\begin{enumerate}
    \item \textbf{AvV}: the full binarization with strengthening solved with  CPLEX default settings.
    \item \textbf{AvV with MIR (resp. $\{0,\frac{1}{2}\}$) cuts only}: the \textbf{AvV} formulation where we set the cut generation to aggressive for MIR (resp. $\{0,\frac{1}{2}\}$) cuts and turn off all other cut types for FCT (resp. CMST).\footnote{The choice of MIR cuts for FCT and $\{0,\frac{1}{2}\}$ cuts (see, \cite{CapraraFischetti1996}) for CMST as leading cut family for the two classes of instances is due to preliminary extended experiments isolating the effect of each family of cutting planes in CPLEX.}
    \item \textbf{AvV with formulation cuts only}: the \textbf{AvV} formulation where we add formulation cuts and solve the formulation with all cuts types turned off using the solver.
    \item \textbf{AvV with formulation and other cuts}: the \textbf{AvV} formulation amended by formulation cuts and solved using the default cut generation settings.
\end{enumerate}

The same experiments are performed with Gurobi for FCT (Table \ref{tab:avv_mir_form_gurobi}) and CMST (Table \ref{tab:avv_MIR2_form_gurobi}), the only difference being that MIR cuts are always for Gurobi the leading cut family again based on extensive computational experiments.


The formulation comparison results are reported in the Tables~\ref{tab:avv_mir_form}--
\ref{tab:avv_MIR2_form_gurobi}. Table~\ref{tab:avv_mir_form} contains the following metrics: (1) the instances $(n,\,B, \, s)$; (2) Root gap\%, which is computed as
$100 \times (\text{OPT}-\text{Root value})/\text{OPT}$; (3) the number of branch-and-bound nodes explored; and (4) the solution time in CPU seconds. The LP gap\% is used as a benchmark. Tables~\ref{tab:avv_zerohalf_form}--\ref{tab:avv_MIR2_form_gurobi} report the same metrics as Table~\ref{tab:avv_mir_form}. We enforce a time limit of 3,600 CPU seconds for each instance. All experiments were run on a computer with an Intel(R) Xeon(R) Gold 6258R processor running at 2.7 GHz, with up to 32 threads used. All optimization problems were solved using the optimization solvers CPLEX with version 22.1.2 and Gurobi with version 13.0.

\setlength{\tabcolsep}{1.3 pt} 
\renewcommand{\arraystretch}{1} 

\begin{table}[htbp]
\centering
\caption{FCT: AvV, MIR, and formulation cuts comparisons using CPLEX.}
\label{tab:avv_mir_form}
{\scriptsize
\begin{tabular}{l@{\hspace{10pt}}r|
                r@{\hspace{8pt}}rr|
                r@{\hspace{8pt}}rr|
                r@{\hspace{8pt}}rr|
                r@{\hspace{8pt}}rr}
\hline
\up\down $(n,\,B,\,s)$ & \begin{tabular}{@{}c@{}}LP\\ gap\%\end{tabular} &
\multicolumn{3}{c|}{AvV} &
\multicolumn{3}{c|}{\begin{tabular}{@{}c@{}}AvV with \\ MIR cuts only\end{tabular}} &
\multicolumn{3}{c|}{\begin{tabular}{@{}c@{}}AvV with\\ formulation cuts only\end{tabular}} &
\multicolumn{3}{c}{\begin{tabular}{@{}c@{}}AvV with formulation\\ and other cuts\end{tabular}} \\
\cline{3-14}
\up\down & &
\begin{tabular}{@{}c@{}}Root\\ gap\%\end{tabular} & Nodes & Time &
\begin{tabular}{@{}c@{}}Root\\ gap\%\end{tabular} & Nodes & Time &
\begin{tabular}{@{}c@{}}Root\\ gap\%\end{tabular} & Nodes & Time &
\begin{tabular}{@{}c@{}}Root\\ gap\%\end{tabular} & Nodes & Time \\
\hline
(30,\,10,\,1) & 13.73 & 0.00 & 0 & 0.78 & 0.00 & 0 & 0.85 & 0.00 & 0 & 0.48 & 0.00 & 0 & 0.56 \\
(30,\,10,\,2) & 14.35 & 0.00 & 0 & 1.24 & 0.00 & 0 & 1.10 & 0.95 & 2,319 & 1.37 & 0.00 & 0 & 0.68 \\
(30,\,10,\,3) & 15.79 & 0.00 & 0 & 1.27 & 0.00 & 0 & 1.03 & 1.27 & 578 & 0.71 & 0.00 & 0 & 0.87 \\
(30,\,10,\,4) & 12.35 & 0.00 & 0 & 0.73 & 1.09 & 58 & 0.80 & 1.07 & 101 & 0.69 & 0.00 & 0 & 0.59 \\
(30,\,10,\,5) & 12.61 & 0.00 & 0 & 0.71 & 0.00 & 0 & 0.80 & 0.35 & 11 & 0.86 & 0.00 & 0 & 0.39 \\
\hline
(30,\,20,\,1) & 15.77 & 0.00 & 0 & 3.75 & 1.60 & 735 & 4.62 & 2.68 & 10,521 & 5.06 & 0.00 & 0 & 2.67 \\
(30,\,20,\,2) & 13.41 & 0.00 & 0 & 1.03 & 0.00 & 0 & 1.67 & 1.17 & 261 & 1.43 & 0.00 & 0 & 1.38 \\
(30,\,20,\,3) & 14.04 & 0.00 & 0 & 1.26 & 0.00 & 0 & 1.35 & 0.69 & 100 & 1.54 & 0.00 & 0 & 0.92 \\
(30,\,20,\,4) & 13.51 & 0.00 & 0 & 2.73 & 1.55 & 246 & 3.37 & 2.44 & 9,859 & 4.66 & 0.00 & 0 & 2.27 \\
(30,\,20,\,5) & 14.72 & 0.00 & 0 & 2.38 & 0.31 & 35 & 2.20 & 1.12 & 598 & 1.52 & 0.00 & 0 & 1.55 \\
\hline
(40,\,10,\,1) & 12.62 & 0.00 & 0 & 1.57 & 0.40 & 42 & 2.33 & 0.49 & 210 & 4.45 & 0.00 & 0 & 0.91 \\
(40,\,10,\,2) & 14.19 & 0.00 & 0 & 1.91 & 1.23 & 1,031 & 4.18 & 2.02 & 10,852 & 4.45 & 0.01 & 9 & 1.75 \\
(40,\,10,\,3) & 11.63 & 0.00 & 0 & 1.31 & 0.00 & 0 & 2.03 & 0.81 & 151 & 1.10 & 0.00 & 0 & 0.86 \\
(40,\,10,\,4) & 10.32 & 0.00 & 0 & 1.45 & 0.00 & 0 & 0.99 & 0.00 & 0 & 1.48 & 0.00 & 0 & 0.57 \\
(40,\,10,\,5) & 11.22 & 0.01 & 0 & 1.76 & 0.00 & 0 & 2.12 & 0.96 & 363 & 2.23 & 0.00 & 0 & 0.82 \\
\hline
(40,\,20,\,1) & 14.62 & 0.00 & 0 & 4.49 & 0.70 & 532 & 7.66 & 1.19 & 9,087 & 7.65 & 0.00 & 0 & 3.43 \\
(40,\,20,\,2) & 16.59 & 0.11 & 7 & 6.60 & 0.49 & 216 & 5.18 & 1.51 & 14,197 & 15.74 & 0.10 & 9 & 5.80 \\
(40,\,20,\,3) & 16.45 & 0.00 & 0 & 7.58 & 1.34 & 1,176 & 9.12 & 2.83 & 75,106 & 43.61 & 0.00 & 0 & 5.52 \\
(40,\,20,\,4) & 12.04 & 0.01 & 0 & 4.57 & 0.92 & 275 & 3.85 & 1.48 & 3,052 & 4.77 & 0.29 & 93 & 3.63 \\
(40,\,20,\,5) & 15.35 & 0.18 & 31 & 9.19 & 0.35 & 63 & 8.26 & 1.43 & 8,974 & 9.25 & 0.00 & 0 & 3.78 \\
\hline
\hline
\textbf{Average} & 13.77 & 0.02 & 1.9 & 2.82 & 0.50 & 220.5 & 3.18 & 1.22 & 7,317.0 & 5.65 & 0.02 & 5.1 & 1.95 \\
\hline
\end{tabular}
}\\[2pt]
\end{table}

\setlength{\tabcolsep}{1.3 pt} 
\renewcommand{\arraystretch}{1} 


\setlength{\tabcolsep}{1.3 pt} 
\renewcommand{\arraystretch}{1} 

\begin{table}[htbp]
\centering
\caption{CMST: AvV, $\{0,\frac{1}{2}\}$, and formulation cuts comparisons using CPLEX.}
\label{tab:avv_zerohalf_form}
{\scriptsize
\begin{tabular}{l@{\hspace{10pt}}r|
                r@{\hspace{8pt}}rr|
                r@{\hspace{8pt}}rr|
                r@{\hspace{8pt}}rr|
                r@{\hspace{8pt}}rr}
\hline
\up\down Instance & \begin{tabular}{@{}c@{}}LP\\ gap\%\end{tabular} &
\multicolumn{3}{c|}{AvV} &
\multicolumn{3}{c|}{\begin{tabular}{@{}c@{}}AvV with \\ $\{0,\frac{1}{2}\}$ cuts only\end{tabular}} &
\multicolumn{3}{c|}{\begin{tabular}{@{}c@{}}AvV with\\ formulation cuts only\end{tabular}} &
\multicolumn{3}{c}{\begin{tabular}{@{}c@{}}AvV with formulation\\ and other cuts\end{tabular}} \\
\cline{3-14}
\up\down & &
\begin{tabular}{@{}c@{}}Root\\ gap\%\end{tabular} & Nodes & Time &
\begin{tabular}{@{}c@{}}Root\\ gap\%\end{tabular} & Nodes & Time &
\begin{tabular}{@{}c@{}}Root\\ gap\%\end{tabular} & Nodes & Time &
\begin{tabular}{@{}c@{}}Root\\ gap\%\end{tabular} & Nodes & Time \\
\hline
tc80-1 & 6.82 & 1.34 & 8,341 & 18.83 & 1.91 & 8,432 & 19.85 & 0.52 & 12,457 & 7.53 & 0.31 & 1,263 & 5.90 \\
tc80-2 & 5.78 & 1.15 & 6,336 & 22.81 & 1.69 & 7,701 & 19.72 & 0.52 & 4,717 & 4.40 & 0.29 & 443 & 6.03 \\
tc80-3 & 5.82 & 0.85 & 8,314 & 13.16 & 1.19 & 6,996 & 16.73 & 0.75 & 8,272 & 3.92 & 0.00 & 0 & 2.11 \\
tc80-4 & 5.56 & 0.37 & 809   & 11.97 & 1.54 & 8,950 & 16.72 & 0.41 & 1,138 & 1.79 & 0.25 & 238 & 3.72 \\
tc80-5 & 6.45 & 0.77 & 5,670 & 13.60 & 1.19 & 10,793& 18.62 & 0.45 & 2,198 & 5.60 & 0.00 & 0 & 5.75 \\
\hline
te80-1 & 6.01 & 1.04 & 5,882 & 17.68 & 1.70 & 8,424 & 22.88 & 1.70 & 194,545 & 38.83 & 0.68 & 10,954 & 13.88 \\
te80-2 & 6.92 & 1.15 & 8,887 & 21.66 & 1.49 & 7,799 & 17.03 & 2.06 & 197,947 & 44.36 & 1.11 & 8,303 & 15.46 \\
te80-3 & 6.13 & 1.26 & 63,266& 56.00 & 2.07 & 71,627& 52.26 & 1.64 & 1,208,665 & 435.46 & 0.83 & 61,801 & 37.05 \\
te80-4 & 7.35 & 1.72 & 92,253& 59.26 & 2.10 & 126,735& 69.72 & 1.51 & 392,730 & 64.76 & 1.21 & 73,564 & 52.64 \\
te80-5 & 7.11 & 1.59 & 7,267 & 21.87 & 2.71 & 33,787 & 26.93 & 2.40 & 678,046 &155.30 & 1.22 & 26,830 & 27.80 \\
\hline
\hline
\textbf{Average} & 6.40 & 1.12 & 20,702.5 & 25.68 & 1.76 & 29,124.4 & 28.05 & 1.20 & 270,071.5 & 76.19 & 0.59 & 18,339.6 & 17.03 \\
\hline
\end{tabular}
}
\end{table}

\setlength{\tabcolsep}{1.3 pt} 
\renewcommand{\arraystretch}{1} 

The results in Tables~\ref{tab:avv_mir_form}--\ref{tab:avv_MIR2_form_gurobi} demonstrate that separating formulation cuts first and then solving the FCT and CMST problems with CPLEX and Gurobi default settings leads to high-quality performance. In particular, the results in Table~\ref{tab:avv_mir_form} show that for FCT, the performance of solving (AvV) with formulation cuts and other cuts is similar to the performance of solving (AvV) with default settings. As an example, the (AvV) formulation with formulation cuts and other cuts solves the particular instance of (40, 20, 5) more effectively by closing the gap to 0.00\% and therefore solving the instance at the root node. The results in Table~\ref{tab:avv_mir_form_gurobi} show that adding formulation cuts to (AvV) and solving the formulation with Gurobi default settings leads to even better performance than solving (AvV) with default settings, as evidenced by the significant difference in root gaps, numbers of explored nodes, and computing times. For the CMST problem, the results in Table~\ref{tab:avv_zerohalf_form} show that CPLEX solves the (AvV) formulation with formulation cuts separated first and other cuts separated later significantly more effectively than the (AvV) formulation with default settings, quantified by the difference between the root gaps of 0.59\% and 1.12\%, numbers of explored nodes, and the computing times. Table~\ref{tab:avv_MIR2_form_gurobi} reports the similarity of performance using Gurobi between the (AvV) formulation with default settings and with formulation cuts and other cuts. These results imply that using the family of formulation cuts at the beginning of the cutting plane separation process can lead to a performance that is in general better than solely using the default cut separation routines of modern MIP solvers such as CPLEX and Gurobi.

\begin{table}[htbp]
\centering
\caption{FCT: AvV, MIR, and formulation cuts comparisons using Gurobi.}
\label{tab:avv_mir_form_gurobi}
{\scriptsize
\begin{tabular}{l@{\hspace{10pt}}r|
                r@{\hspace{8pt}}rr|
                r@{\hspace{8pt}}rr|
                r@{\hspace{8pt}}rr|
                r@{\hspace{8pt}}rr}
\hline
\up\down $(n,\,B,\,s)$ & \begin{tabular}{@{}c@{}}LP\\ gap\%\end{tabular} &
\multicolumn{3}{c|}{AvV} &
\multicolumn{3}{c|}{\begin{tabular}{@{}c@{}}AvV with \\ MIR cuts only\end{tabular}} &
\multicolumn{3}{c|}{\begin{tabular}{@{}c@{}}AvV with\\ formulation cuts only\end{tabular}} &
\multicolumn{3}{c}{\begin{tabular}{@{}c@{}}AvV with formulation\\ and other cuts\end{tabular}} \\
\cline{3-14}
\up\down & &
\begin{tabular}{@{}c@{}}Root\\ gap\%\end{tabular} & Nodes & Time &
\begin{tabular}{@{}c@{}}Root\\ gap\%\end{tabular} & Nodes & Time &
\begin{tabular}{@{}c@{}}Root\\ gap\%\end{tabular} & Nodes & Time &
\begin{tabular}{@{}c@{}}Root\\ gap\%\end{tabular} & Nodes & Time \\
\hline
(30,\,10,\,1) & 12.96 & 0.00 & 0 & 1.46 & 0.00 & 0 & 0.94 & 0.00 & 0 & 0.46 & 0.00 & 0 & 1.08 \\
(30,\,10,\,2) & 14.16 & 0.00 & 0 & 2.86 & 0.00 & 0 & 1.34 & 0.00 & 0 & 0.50 & 0.00 & 0 & 0.69 \\
(30,\,10,\,3) & 15.29 & 0.00 & 0 & 2.42 & 0.00 & 0 & 1.39 & 1.16 & 288 & 0.87 & 0.00 & 0 & 0.96 \\
(30,\,10,\,4) & 11.95 & 0.00 & 0 & 1.56 & 0.00 & 0 & 1.25 & 0.00 & 0 & 0.46 & 0.00 & 0 & 0.65 \\
(30,\,10,\,5) & 12.54 & 0.00 & 0 & 1.19 & 0.00 & 0 & 0.99 & 0.00 & 0 & 0.91 & 0.00 & 0 & 0.72 \\
\hline
(30,\,20,\,1) & 15.63 & 0.13 & 10 & 6.62 & 0.60 & 166 & 4.12 & 2.68 & 8,235 & 4.43 & 0.00 & 0 & 3.10 \\
(30,\,20,\,2) & 13.31 & 0.00 & 0 & 2.90 & 0.00 & 0 & 2.01 & 0.83 & 215 & 1.75 & 0.00 & 0 & 1.81 \\
(30,\,20,\,3) & 13.86 & 0.00 & 0 & 3.29 & 0.00 & 0 & 1.80 & 0.00 & 0 & 0.91 & 0.00 & 0 & 1.98 \\
(30,\,20,\,4) & 13.36 & 0.00 & 0 & 4.87 & 0.71 & 180 & 3.84 & 2.48 & 5,281 & 3.02 & 0.00 & 0 & 2.83 \\
(30,\,20,\,5) & 14.48 & 0.00 & 0 & 3.87 & 0.00 & 0 & 1.79 & 1.12 & 1,201 & 1.49 & 0.00 & 0 & 2.79 \\
\hline
(40,\,10,\,1) & 12.42 & 0.00 & 0 & 2.65 & 0.00 & 0 & 1.60 & 0.31 & 42 & 1.04 & 0.00 & 0 & 1.19 \\
(40,\,10,\,2) & 13.90 & 0.00 & 0 & 4.30 & 0.69 & 152 & 3.12 & 1.60 & 3,816 & 2.57 & 0.00 & 0 & 2.39 \\
(40,\,10,\,3) & 10.90 & 0.00 & 0 & 3.39 & 0.00 & 0 & 1.68 & 0.00 & 0 & 0.96 & 0.00 & 0 & 1.47 \\
(40,\,10,\,4) & 9.53 & 0.00 & 0& 2.04 & 0.00 & 0 & 1.88 & 0.00 & 0 & 0.72 & 0.00 & 0 & 1.40 \\
(40,\,10,\,5) & 11.13 & 0.00 & 0 & 2.08 & 0.00 & 0 & 1.65 & 0.88 & 129 & 0.97 & 0.00 & 0 & 1.19 \\
\hline
(40,\,20,\,1) & 14.58 & 0.01 & 0 & 13.40 & 0.60 & 192 & 7.35 & 1.11 & 7,706 & 4.43 & 0.00 & 0 & 5.22 \\
(40,\,20,\,2) & 16.46 & 0.08 & 4 & 9.90 & 0.29 & 97 & 10.59 & 1.41 & 8,588 & 11.01 & 0.01 & 0 & 5.17 \\
(40,\,20,\,3) & 16.24 & 0.64 & 621 & 11.26 & 0.84 & 635& 6.92 & 2.82 & 34,645& 29.74 & 0.46 & 194 & 7.13 \\
(40,\,20,\,4) & 11.71 & 0.38 & 30 & 13.45 & 0.82 & 350 & 9.42 & 1.14 & 1,596 & 2.87 & 0.00 & 0 & 6.51 \\
(40,\,20,\,5) & 15.34 & 0.00 & 0 & 16.94 & 0.00 & 0& 7.88 & 0.98 & 5,631 & 4.89 & 0.01 & 0 & 6.60 \\
\hline
\hline
\textbf{Average} & 13.49 & 0.06 & 33.25 & 5.52 & 0.23 & 88.6 & 3.58 & 0.93 & 3,868.7 & 3.70 & 0.02 & 9.7 & 2.74 \\
\hline
\end{tabular}
}\\[2pt]
\end{table}
As reported in Tables~\ref{tab:avv_mir_form} and~\ref{tab:avv_mir_form_gurobi}, for the FCT problem, MIR cuts close the gap to the average value of 0.50\% and 0.23\% using CPLEX and Gurobi, respectively. In Table~\ref{tab:avv_mir_form}, we observe that it is significantly more effective to solve instances with MIR cuts only than with formulation cuts only in CPLEX, as highlighted by the difference between the average root gaps of 0.50\% and 1.22\% and the average number of explored nodes of 220.5 and 7,317.0, respectively. Similarly, using Gurobi, it is significantly more effective to solve instances with MIR cuts only than with formulation cuts only, as shown by the difference between the average root gaps of 0.23\% and 0.93\% and the average number of explored nodes of 88.6 and 3,868.7 in Table~\ref{tab:avv_mir_form_gurobi}. 

\begin{table}[htbp]
\centering
\caption{CMST: AvV, MIR, and formulation cuts comparisons using Gurobi.}
\label{tab:avv_MIR2_form_gurobi}
{\scriptsize
\begin{tabular}{l@{\hspace{10pt}}r|
                r@{\hspace{8pt}}rr|
                r@{\hspace{8pt}}rr|
                r@{\hspace{8pt}}rr|
                r@{\hspace{8pt}}rr}
\hline
\up\down Instance & \begin{tabular}{@{}c@{}}LP\\ gap\%\end{tabular} &
\multicolumn{3}{c|}{AvV} &
\multicolumn{3}{c|}{\begin{tabular}{@{}c@{}}AvV with \\ MIR cuts only\end{tabular}} &
\multicolumn{3}{c|}{\begin{tabular}{@{}c@{}}AvV with\\ formulation cuts only\end{tabular}} &
\multicolumn{3}{c}{\begin{tabular}{@{}c@{}}AvV with formulation\\ and other cuts\end{tabular}} \\
\cline{3-14}
\up\down & &
\begin{tabular}{@{}c@{}}Root\\ gap\%\end{tabular} & Nodes & Time &
\begin{tabular}{@{}c@{}}Root\\ gap\%\end{tabular} & Nodes & Time &
\begin{tabular}{@{}c@{}}Root\\ gap\%\end{tabular} & Nodes & Time &
\begin{tabular}{@{}c@{}}Root\\ gap\%\end{tabular} & Nodes & Time \\
\hline
tc80-1 & 6.82 & 0.38 & 2,993 & 5.95 & 0.69 & 7,862 & 5.12 & 0.38 & 4,622 & 8.31 & 0.38 & 1,428 & 6.75 \\
tc80-2 & 5.77 & 0.37 & 4,285 & 7.67 & 0.67 & 6,274 & 5.27 & 0.30 & 2,000 & 2.33 & 0.30 & 4,142 & 6.71 \\
tc80-3 & 5.79 & 0.21 & 52 & 3.86 & 0.41 & 962 & 3.95 & 0.69 & 1,352 & 1.84 & 0.21 & 169 & 3.61 \\
tc80-4 & 5.53 & 0.32 & 691 & 3.35 & 0.40 & 2,787 & 3.91 & 0.32 & 570 & 1.91 & 0.00 & 0 & 3.43 \\
tc80-5 & 6.41 & 0.18 & 80 & 4.84 & 0.36 & 277 & 3.90 & 0.30 & 3,624 & 3.97 & 0.00 & 0 & 4.82 \\
\hline
te80-1 & 5.97 & 1.05 & 11,790 & 6.30 & 1.31 & 16.294 & 7.21 & 1.67 & 88,847 & 19.57 & 0.58 & 3,929 & 8.05 \\
te80-2 & 6.91 & 1.10 & 13,544 & 8.51 & 1.20 & 11,147 & 7.94 & 1.31 & 143,357 & 24.43 & 0.68 & 3,529 & 11.57 \\
te80-3 & 6.10 & 1.00 & 20,348 & 28.91 & 1.43 & 88.543 & 44.45 & 1.59 & 810,164 & 125.23 & 0.64 & 14,497& 22.08 \\
te80-4 & 7.33 & 1.18 & 78,671 & 28.74 & 1.55 & 329,748 & 53.32 & 1.50 & 361,875 & 55.87 & 1.05 & 66,907& 46.88 \\
te80-5 & 7.08 & 1.08 & 6,402& 14.81 & 1.33 & 27,002 & 14.85 & 2.36 & 417,130 & 73.77 & 1.23 & 15,636 & 36.15 \\
\hline
\hline
\textbf{Average} & 6.37 & 0.69 & 13,885.6 & 11.29 & 0.94 & 49,089.6 & 14.98 & 1.04 & 183,354.1 & 31.72 & 0.51 & 11,023.7 & 15.00 \\
\hline

\end{tabular}
}\\[2pt]
\end{table}

For the CMST problem, Tables~\ref{tab:avv_zerohalf_form} and~\ref{tab:avv_MIR2_form_gurobi}, $\{0,\frac{1}{2}\}$ cuts contribute significantly to closing the gap to the average value of 1.76\% using CPLEX, and MIR cuts contribute as the leading family to closing the gap to the average value of 0.94\% using Gurobi, respectively. In CPLEX, one gets a significantly lower number of explored nodes, lower computing time, and a slightly higher root gap\% using $\{0,\frac{1}{2}\}$ cuts only instead of formulation cuts only, as demonstrated in Table~\ref{tab:avv_zerohalf_form}. For Gurobi, Table~\ref{tab:avv_MIR2_form_gurobi} suggests that using MIR cuts only instead of formulation cuts only is more effective on average.

\section{Conclusions}
\label{sec:conclusions}

In this paper, we have shown that the use of extended formulations replacing bounded general integer variables by linear combinations of a set of auxiliary binary variables linked by additional linear constraints can be highly beneficial for enhancing the performance of MIP solvers. A similar message has been occasionally reported by several authors in the past (for example, \cite{U2011}, \cite{BM2015}, \cite{AVV2017}) lacking of an extensive computational analysis about which, when and why extended formulations could make a difference. 

We concentrated on MIP instances belonging to the general class of network flow problems and we have shown that the appeal of some of those extended formulations is in triggering effective cut generation by MIP solvers. Moreover, we have shown that a specific simple family of mixed-integer rounding inequalities known as formulation cuts is especially beneficial if separated on a specific group of constraints of the network flow problems, flow-conservation inequalities. We have computationally pinpointed that these constraints reformulated in terms of binary variables trigger cut generation in MIP solvers by analyzing extensively the results for both CPLEX and Gurobi. 

To the best of our knowledge, this is the first extensive computational analysis of the topic and the results significantly contribute to clarify when and how binarizations can be effective in MIP with network flow structure, thus giving practical guidelines for their use.


%

%
%

\ACKNOWLEDGMENT{We would like to express our sincere gratitude to Andrea Tramontani and Oktay G\"unl\"uk for extensive discussions on the topic in early stages of this research. Thanks are also due to Mathieu Van Vyve and Gustavo Angulo for providing the generator of FCT instances.}







\bibliographystyle{informs2014}
\bibliography{literature}

\end{document}